\documentclass[11pt]{amsart}
\usepackage{amssymb}
 \newcommand{\Z}{{\mathbb Z}}

 \newcommand{\Q}{{\mathbb Q}}
 \newcommand{\R}{{\mathbb R}}
 \newcommand{\K}{{\mathbb K}}

 \newtheorem{theorem}{Theorem}[section]
 \newtheorem{lemma}[theorem]{Lemma} 
 
 \newtheorem{proposition}[theorem]{Proposition}
 \newtheorem{corollary}[theorem]{Corollary}

 \newtheorem{definition}[theorem]{Definition}   
 \newtheorem{conjecture}[theorem]{Conjecture}

 \def\Box
  {\hfill \thinspace\vbox{\hrule height .5pt \hbox{\vrule  
   width .5pt \vbox to 7pt{\hbox to 3.5pt{}} \vrule width .5pt} 
   \hrule height 0pt depth .5pt}}

\title{Topology on the spaces of orderings of groups}

\author{Adam S. Sikora}

\date{}
\begin{document}
\begin{abstract}
A natural topology on the space of
left orderings of an arbitrary semi-group is introduced.
It is proved that this space is compact and that for free 
abelian groups it is homeomorphic to the Cantor set.
An application of this result is a new proof of the existence of 
universal Gr\"obner bases. 

\end{abstract} 

\maketitle
\section{Orderings for semi-groups}
Given a semi-group $G$ (ie. a set with an associative binary operation),
a linear order, $<,$ on $G$ is a left order if $a< b$ implies
$ca< cb,$ for any $c.$ Similarly, a linear order $,<,$
is a right order if $a< b$ implies $ac< bc,$ for any $c\in G.$
The sets of all left and right orderings of $G$ are denoted by $LO(G)$ and
$RO(G)$ respectively.
If $G$ is a group then there is a 1-1 correspondence between these two sets
which associates with any left ordering, $<_l,$ a right
ordering, $<_r,$ such that $a<_r b$ if and only if $b^{-1}<_l a^{-1}.$
For more about ordering of groups see \cite{Fu,KK,MR}.

Let $U_{a,b}\subset LO(G)$ denote the set of all left orderings, $<,$ for
which $a< b.$ We can put a topology on $LO(G)$ in one of the following two 
ways.

\begin{definition}\label{d1} $LO(G)$ has the
smallest topology for which all the sets $U_{a,b}$ are open. Any open set 
in this topology is a union of sets of the form $U_{a_1,b_1}\cap ...\cap 
U_{a_n,b_n}.$
\end{definition}

\begin{definition}\label{d2} Let $G_0\subset G_1\subset G_2 ...\subset G$
be an arbitrary complete filtration of $G$ by its subsets.  
(A filtration is complete if $\bigcup_i G_i=G$). For $<_1,<_2\in LO(G)$ 
we define $\rho(<_1,<_2)$ to be
$\frac{1}{2^r},$ where $r$ is the largest
number with the property that $<_1$ and $<_2$ coincide when restricted
to $G_r;$ We put \mbox{$\rho(<_1,<_2)=0$} if such $r$ 
does not exist ($r=\infty$).
\end{definition}

From now on we will consider countable semi-groups $G$ only and
such filtrations only which are composed of finite subsets of $G.$

\begin{proposition}\label{metric=topology} $\rho$ is a metric on $LO(G)$ and
the topology on $LO(G)$ induced by that metric
coincides with the topology introduced in Definition \ref{d1}. 
In particular, it does not depend on the choice of a filtration of $G.$
\end{proposition}

\begin{proof} It is easy to check that $\rho$ is a metric.
Hence, the proposition follows from the following two statements:\\
(1) any open ball $B(<_0,1/2^r)$ (with respect to the metric
$\rho$) is open in the topology introduced in Definition \ref{d1};\\
{\it Proof:} $<_1\in B(<_0,1/2^r)$ if and only if $<_1$ and 
$<_0$
coincide on the set $G_{r+1}.$ Therefore $$B(<_0,1/2^r)=\bigcap U_{ab},$$
where the intersection is taken over all possible pairs $(a,b)$ of elements 
of $G_{r+1}$ for which $a<_0 b.$\\
(2) any set of the form $U_{a_1b_1}\cap U_{a_2b_2}\cap ... \cap U_{a_nb_n}$
is open with respect to the metric $\rho.$ In other words,
for any $<_0\in U_{a_1b_1}\cap U_{a_2b_2}\cap ... \cap U_{a_nb_n}$ there
exist $r$ such that $B(<_0,1/2^r)\subset U_{a_1b_1}\cap U_{a_2b_2}\cap ... 
\cap U_{a_nb_n}.$\\
{\it Proof:} There exists an element of the filtration, $G_r,$ such 
that\\
$a_1,b_1,...,a_n,b_n\in G_r.$ For such $r,$ $B(<_0,1/2^r)\subset U_{a_1b_1}
\cap U_{a_2b_2}\cap ... \cap U_{a_nb_n}.$
\end{proof}

Recall that a space is totally disconnected if every 
two distinct points of it are contained in two disjoint open sets 
covering the space.

\begin{theorem}\label{compact} $LO(G)$ is a compact, totally disconnected 
topological space. 
\end{theorem}

\noindent {\it Proof:}
For any two left orderings $<_1,<_2\in LO(G)$ there exist $a,b\in G$ 
such that
$<_1\in U_{ab}$ and $<_2\in U_{ba}.$ Since $U_{ab}\cup U_{ba}=LO(G),$
$U_{ab}\cap U_{ba}=\emptyset,$ $LO(G)$ is totally disconnected.
Now we are going to show that $LO(G)$ is compact.
Consider any complete, infinite filtration of $G$ by its finite subsets
and the associated metric $\rho.$ We need to prove that any infinite 
sequence $<_1,<_2,...\in LO(G)$ has a convergent subsequence. 
We construct this subsequence in the following manner: 
Since there are only finitely many possible linear orderings of elements 
of $G_1$ there is an infinite subsequence $<_{i_1^1},<_{i_2^1},
<_{i_3^1}...$ of $<_1,<_2,...$ whose elements induce the same 
linear order on $G_1.$
Now, pick out of this sequence an infinite subsequence of orders,
$<_{i_1^2},<_{i_2^2},...,$ which agree on
$G_2.$ Continue this process for $G_3,G_4,...$ {\it ad infinitum.}
Consider now a sequence $<^1,<^2,<^3,....$ constructed by
picking up the $n$-th element from the $n$-th subsequence constructed 
above for $n=1,2,...$ Since $<^1,<^2,<^3,....$ is a subsequence of 
$<_1,<_2,...$ the following lemma completes the proof.

\begin{lemma} $<^1,<^2,<^3,....$ converges to a left order $<^{\infty},$ 
defined as follows: $a<^{\infty} b$ if and only if $a<^n b$ for 
almost all $n.$
\end{lemma}

\begin{proof}
If $a,b\in G_r$ then either $a<^i b$ for $i>r$ or $b<^i a$ for $i>r.$
Therefore, $<^{\infty}$ is a total order and, it is easy to verify that
it is also a left order on $G.$ Since 
$\rho(<^n,<^{\infty})\leq \frac{1}{2^n},$ the sequence 
$<^1,<^2,<^3,....$ converges to $<^{\infty}.$
\end{proof}

\begin{corollary}\label{cor_for_compact}
$LO(G)$ is homeomorphic to the Cantor set if and only if\\
(1) $LO(G)\ne \emptyset,$ and\\
(2) for any sequence $a_1,b_1,....,a_n,b_n,$ of elements of $G,$ the set
$U_{a_1b_1}\cap U_{a_2b_2}\cap ... \cap U_{a_nb_n}$ is either empty or 
infinite.
\end{corollary}

\begin{proof}
Any nonempty, metrizable, compact, perfect 
 and totally disconnected set is a 
Cantor set (\cite[Corollary 2-98]{HY}). 
(A set is perfect if every point of the set is a limit point.)
Condition (2) states that $LO(G)$
is perfect.
\end{proof}

Orderings on $\Z^n$ were classified in \cite{Te, R}.

\begin{proposition}\label{Z^n}
For $n>1,$ $LO(\Z^n)$ is homeomorphic to the Cantor
set.
\end{proposition}

\begin{proof} 
Consider the smallest $n>1$ for
which the statement fails. By Corollary 
\ref{cor_for_compact}, there exists a finite set of pairs
$(a_1,b_1),...,(a_s,b_s)\in \Z^n\times \Z^n$ such that
the number of orderings $<$ on $\Z^n$ such that
$a_i<b_i,$ for $i=1,...,s,$ is positive and finite.
By adding some additional pairs, if necessary,
we can assume that there is only one such ordering, $<.$
Furthermore, we can assume that $b_j-a_j$ is not a rational multiple
of $b_i-a_i,$ for any $i\ne j.$
We extend $<$ to
an ordering on $\Q^n$ by demanding that $v_1<v_2$ for $v_1,v_2\in \Q^n$ if
and only if $nv_1<nv_2$ for all $n\in \Z$ such that $nv_1,nv_2\in \Z^n.$
Consider the set $H\subset \Q^n\otimes \R=\R^n$ composed of elements
$x$ such that any neighborhood of $x$ in $\Q^n$ (with respect to the
Euclidean topology in $\R^n$) contains both
positive and negative elements. ($v\in \Q^n$ is positive if $0<v$).
One can prove that $H$ is a hyper-plane
in $\R^n$ and that the two connected components of $\R^n\setminus H,$ 
denoted by $H_+$ and $H_-,$ have the property that all elements of $\Q^n$ in
$H_+$ are positive and all elements of $\Q^n$ in $H_-$ are negative.
Therefore the vectors $b_i-a_i$ lie either in $H_+$ or in $H.$
Denote by $I$ the set of $i$'s such that $b_i-a_i\in H.$
Observe that $<$ is the only order on $H\cap \Z^n$ for which 
$a_i<b_i$ for $i\in I.$
(Any other order on $H\cap \Z^n$ satisfying this condition would extended to
an order $<'$ on $\Q^n$ for which $a_i<'b_i$ for all $i=1,...,k.$ This would
contradict the uniqueness of $<$). 
Since $H\cap \Z^n$ does not satisfy condition (2) of Corollary 
\ref{cor_for_compact}, the initial assumption about $n$ implies that
$H\cap \Z^n= \emptyset$ or $\Z.$ If $H\cap \Z^n= \emptyset$ then
$I=\emptyset$ and there exist infinitely many hyper-planes 
$H'\subset \R^n,$ obtained by small perturbations of $H$, giving
rise to infinitely many orderings $<'$ of $\Q^n$ such that 
$a_i<' b_i,$ for $i=1,...,k.$
Therefore, $H\cap \Z^n=\Z.$ Since we assumed that none of the
vectors $b_i-a_i$ is a rational multiple of another, we see 
that $b_i-a_i$ lies in $H$ for exactly one index $i=i_0.$ Now, by small
perturbations of $H$ in $\R^n$ we can obtain infinitely many new hyper-planes
$H'$ such that $b_{i_0}-a_{i_0}$ lies in the same component
of $\R^n\setminus H'$ as all other vectors $b_i-a_i.$
Each of these hyper-planes induces an ordering $<'$ on $\Z^n$ such that
$a_i<'b_i$ for $i=1,...,k$ -- a contradiction.
\end{proof}

\section{Bi-orderings}
A bi-ordering is a linear ordering which is both left and right ordering.
Using Proposition \ref{metric=topology}, we see that the set of
bi-orderings $BiO(G)$ inherits the same topology from $LO(G)$ and
$RO(G).$ 

\begin{proposition}\label{bio} $BiO(G)$ is a closed subset of $LO(G).$
Hence $BiO(G)$ is the Cantor set if and only if
$BiO(G)\ne \emptyset$ and for any sequence $a_1,b_1,....,a_n,b_n,$ 
of elements of $G,$ the set
$BiO(G)\cap U_{a_1b_1}\cap U_{a_2b_2}\cap ... \cap U_{a_nb_n}$
is either empty or infinite.
\end{proposition}

\begin{proof} Consider a metric on $LO(G)$ induced by a filtration of $G.$
For any infinite sequence, $<_1,<_2,...,$ of
bi-orders of $G$ converging to a left order, $<_{\infty},$ we have 
$a<_{\infty} b$ if
and only if $a<_i b$ for almost all $i.$ Hence $a<_{\infty} b$ implies
$ac<_{\infty} bc$ for any $c.$ This proves that $BiO(G)$ is a closed subset 
of $LO(G).$ The second part of the statement is proved exactly as
Corollary \ref{cor_for_compact}.
\end{proof}

One should not expect, however, that the set of bi-orderings of 
a group, if infinite, is homeomorphic to the Cantor set. There are 
known examples of groups $G$ for which $BiO(G)$ has countably 
infinite number of elements, \cite{MR}.

\begin{conjecture}
For the free group on $n>1$ generators, $F_n,$ 
the spaces $LO(F_n)$ and $BiO(F_n)$ are homeomorphic to the Cantor set.
\end{conjecture}

Let $G=G^0\supset G^1\supset ...$ be the lower central series
of a group $G$ such that $\bigcap_{k=0}^\infty G^k=1.$ 
\v Simbireva and Neumann (see \cite{Si, Ne} and \cite[Ch. IV \S 2]{Fu}) 
showed that any choice of orders on the groups 
$G^k/G^{k+1}$ yields a total bi-order on $G:$ 
if $g\in G^k\setminus G^{k+1}$ then $g$ is positive in $G$
if and only if $g$ is positive in $G^{k}/G^{k+1}.$ Since any torsion-free 
abelian group is orderable, \cite{Te}, if the groups $G^k/G^{k+1}$ are 
torsion-free then $G$ is bi-orderable
\footnote{In fact the result of \v Simbireva and Neumann 
applies to all groups with a central series, including the transfinite ones,
ending with the trivial group.}.
We will call such orders  standard and denote their set by
$SBiO(G).$ They are characterized by the following condition satisfied by
all $g\in G$: if $g\in G^k\setminus G^{k+1}$ and $g$ is positive 
then all elements of $gG^{k+1}$ are positive. 

\begin{proposition}\label{standard}
(1) $SBiO(G)$ is a closed subset of $BiO(G).$\\
(2) If $G\ne \Z$ and each factor $G^k/G^{k+1}$ is finitely generated then 
$SBiO(G)$ is either empty or homeomorphic to the Cantor set.
\end{proposition}

It is possible that the results of \cite{Te} make possible to 
relax the assumption about the finite number of generators of
each of the factors $G^k/G^{k+1}.$

\noindent{\it Proof:} (1) If an order $\leq$ is not standard then,
by the above characterization of standard orders, there exists
$g\in G^k\setminus G^{k+1}$ and $h\in G^{k+1}$ such that $g>e$ and
$gh<e.$ Therefore $\leq$ has an open neighborhood composed of 
non-standard orders.

(2) Assume that $G\ne \Z.$ We may also assume that 
all $G^k/G^{k+1}$ are torsion-free and, hence, free abelian groups, since
otherwise $SBiO(G)=\emptyset.$ As before, it suffices to prove that
for any $a_1,b_1,...,a_n,b_n\in G$ the set 
$$SBiO(G)\cap U_{a_1b_1}\cap U_{a_2b_2}\cap ... \cap U_{a_nb_n}$$
is either empty or infinite. Assume that this set is not empty.
If $G^k/G^{k+1}=\Z^n$ for $n\geq 2$ for some $k$ then one can obtain 
infinitely many standard orders in 
$U_{a_1b_1}\cap U_{a_2b_2}\cap ... \cap 
U_{a_nb_n}$ by slightly modifying the order on $G^k/G^{k+1}.$
(We use here the classification of orders of $\Z^n$ given 
in the proof of Proposition \ref{Z^n}). Therefore from now on it is enough
to assume that $G^k/G^{k+1}$ is either $\Z$ or trivial for all $k.$ 

We claim that in this situation $G^k/G^{k+1}= \Z$ for all $k.$ 
Indeed if $G^k/G^{k+1}$ is trivial for some $k$ then 
$G^l=G^k$ for all $l\geq k,$
and since $\bigcap_{i=0}^\infty G^i=1,$ 
$G$ is nilpotent. However, the only nilpotent group $G$
with $G/G^1=\Z$ is $G=\Z.$ (This statement follows from
\cite[Cor. 10.3.3]{Ha} where $H=<h>$ and $hG^1$ is a generator
of $G/G^1=\Z$).

Therefore, we may assume that $G^k/G^{k+1}=\Z$ for all $k.$
Since each $G^k/G^{k+1}=\Z$ has exactly two orderings it is not hard to
see that $SBiO(G)=\{0,1\}^{\aleph_0}$ as a topological space.
\Box

An example of a non-standard order can be constructed as follows:
Vinogradov proved that if $A,B$ are bi-ordered groups then $A*B$ 
has also a bi-order, \cite{Vi,MR}.
The construction of such bi-order given in the proof of
Theorem 2.3.1 in \cite{MR} has the property that if 
$a_0,a_1\in A,$ $b\in B,$ $a_0>e,$ $a_1<e,$ $b>e$ then 
$a_0ba_1b^{-1}a_1^{-1}<e.$ Therefore that order is non-standard.

\section{Example: $\Z\times \Z$}

Let $\R_{[)}$ denote the real line with the topology with a basis of
open sets composed of intervals $[a,b).$ Analogously, let $\R_{(]}$ be the
real line with the topology whose basis is composed of intervals $(a,b].$
The above topologies carry onto the unit circle in $\R^2$
via the map $z\to e^{2\pi i z},$ $\R_{[)}/\Z\to S^1_{[)}\subset \R^2,$
$\R_{(]}/\Z \to S^1_{(]}\subset \R^2.$ 
We say that a point $(x,y)\in S^1\subset \R^2$ is rational if $x/y\in \Q.$

Let $X$ be the union of these two circles with corresponding irrational 
points identified, $X=(S^1_{[)}\sqcup S^1_{(]})/\sim,$ where 
$(x_1,y_1)\sim (x_2,y_2)$ if $(x_1,y_1)\in S^1_{[)},$ 
$(x_2,y_2)\in S^1_{(]},$
and $(x_1,y_1)=(x_2,y_2)$ is irrational.

\begin{proposition}\label{Z^2}
$BiO(\Z^2)$ is homeomorphic to $X.$
\end{proposition}

\begin{proof}
First, we construct a map $S^1_{[)}\sqcup S^1_{(]}\to BiO(\Z^2)$ as follows:
Associate with $x\in S^1_{[)}$ an order $<_x$ on $\Z^2$ such that
$y\in \Z^2$ is positive if and only if the (oriented) angle 
between the vectors $y$ and $x$ in $\R^2$ is in the interval $(0,\pi].$
Similarly, we associate with $x\in S^1_{(]}$ an order 
$<^x$ on $\Z^2$ such that
$y\in \Z^2$ is positive if and only if the angle 
between $y$ and $x$ is in the interval $[0,\pi).$
Note that these maps descend to a 1-1 map $i$ from $X$ into $BiO(G).$

We claim that $i$ is onto: Recall from the proof of 
Proposition \ref{Z^n} that any ordering on $\Z^2$ defines a 
$1$-dimensional subspace $H$ in $\R^2,$ such that one component of 
$\R^2\setminus H$ is composed
of positive elements and the other of negative elements. Notice
that all orderings on $\Z^2$ inducing the line $H$
are of the form $<_x$ and $<^x,$ where $x$ is one of the two unit vectors
in $H.$

Finally, we claim that $i$ is a continuous map. To verify that claim it is
sufficient to check that $i^{-1}(U_{0a})\subset X$ is open for any 
$a\in \Z^2.$
This condition follows from the fact that the sets
$\{x\in S^1_{[)}: 0<_x a\}$ and $\{x\in S^1_{(]}: 0<^x a\}$ are open
in $ S^1_{[)}$ and in $ S^1_{(]},$ respectively.
\end{proof}

The above proposition has an unexpected application.
Although it is known that any compact set is a continuous image of
a map defined on the Cantor set (see \cite[Thm 3-28]{HY}), usually 
it is difficult to visualize such a map.
However, we can visualize a map from the Cantor set onto $S^1.$
Consider the map $f: X\to S^1$ given by identification of all corresponding
points, $x\in S^1_{[)}$ with $x\in S^1_{(]},$ in 
$X=S^1_{[)}\sqcup S^1_{(]}/\sim.$
Since the topology of $S^1_{(]}$ and of $S^1_{[)}$ is richer than the
Euclidean topology of the circle, $f$ is continuous.
By Propositions \ref{Z^n} and \ref{Z^2}, $X$ is the Cantor set.

\section{Applications to Gr\"obner bases}

Let $\K[x_1,...,x_n]$ be a polynomial ring over a field $\K.$
Monomials in it, $x_1^{i_1}....x_n^{i_n},$ form a monoid (a semi-group
with an identity), isomorphic to $\Z_{\geq 0}^n$ by the isomorphism
which carries $x_1^{i_1}....x_n^{i_n}$ to $(i_1,...,i_n).$ From now on
we will identify these two monoids.

A linear order on a set $G$ is a well-order if 
any subset of $G$ has a smallest element.
For a semi-group $G,$ we denote the set of all left well-orders of $G$
by $LWO(G).$ The elements of $LWO(\Z_{\geq 0}^n)$ are called monomial 
orderings on $\K[x_1,...,x_n].$ As a consequence of Dickson's lemma,
\cite[Cor. 6 Ch.2 \S4]{CLO},
an order $<$ on $\Z_{\geq 0}^n$ is a well-order if and only if
$0$ is the smallest element for $<.$ Therefore
$$LWO(\Z_{\geq 0}^n)=LO(\Z_{\geq 0}^n)\setminus \bigcup U_{a0},$$
where the sum is over all non-zero elements of $\Z_{\geq 0}^n.$
Since $\bigcup U_{a0}$ is open, we get the following.

\begin{corollary}\label{well_orders_closed} 
$LWO(\Z_{\geq 0}^n)$ is a closed subset of 
$LO(\Z_{\geq 0}^n).$ Hence, by Theorem \ref{compact},
$LWO(\Z_{\geq 0}^n)$ is compact.
\end{corollary}

Furthermore, by adopting the reasoning given in the proof of Proposition 
\ref{Z^n}, one can show that both $LWO(\Z_{\geq 0}^n)$ and 
$LO(\Z_{\geq 0}^n)$ are homeomorphic to 
the Cantor set for $n>1.$

Any polynomial $w\in \K[x_1,...,x_n]$ decomposes as $\sum_1^d c_im_i,$ where 
$m_i$'s are monomials, $m_i\ne m_j,$ and $c_i$'s are scalars in $\K,$
$c_i\ne 0.$ For a given monomial ordering,
the leading monomial of $w$, $LM(w),$ is the
largest monomial among $m_1,...,m_d.$
If $LM(w)=m_i,$ then the leading term of $w$ is $LT(w)=c_im_i.$
For a non-zero ideal $I\triangleleft \K[x_1,...,x_n],$ we denote
by $LM(I)$ the ideal in $\K[x_1,...,x_n]$ generated by the leading 
monomials of polynomials in $I.$ Finally, we say that a set of
polynomials in $\{f_1,...,f_d\}\subset I,$ is a {\em Gr\"obner basis} of
$I$ if the leading monomials of these polynomials generate $LM(I).$
Such a basis is very useful for practical computations with $I$ and 
$\K[x_1,...,x_n]/I,$ see \cite{AL,BW,CLO}.
Buchberger's algorithm provides a practical method 
of calculating a Gr\"obner bases, see eg. \cite[Ch. 2 \S7]{CLO}.

\begin{proposition}\label{open}
For any ideal $I\triangleleft \K[x_1,...,x_n]$ and for any set of
polynomials $f_1,...,f_s\in I$ the set of monomial orderings
on $\K[x_1,...,x_n]$ for which $f_1,...,f_s$ is a Gr\"obner basis of $I$
is open.
\end{proposition}

\begin{proof}
For any $f,g\in \K[x_1,...,x_n],$ let
$$S(f,g)=\frac{LCM(LM(f),LM(g))}{LT(f)}f-\frac{LCM(LM(f),LM(g))}{LT(g)}g,$$
where $LCM$ denotes the least common multiple.
By \cite[Theorem 2.\S6.6]{CLO}, a set 
$\{f_1,...,f_s\}$ is a Gr\"obner basis for
$I$ if and only if $S(f_i,f_j)$ is divisible by
$\{f_1,...,f_s\}$ (with remainder $=0$). (The division by
multiple polynomials was defined in \cite[Theorem 2.\S3.3]{CLO}).
The division process 
is based on a finite number of comparisons $m\leq m'$ between monomials.
Therefore, if $S(f_i,f_j)$ is divisible by
$\{f_1,...,f_s\}$ for a monomial ordering $\leq,$ then
it is also divisible by $\{f_1,...,f_s\}$ for all orderings in some
open neighborhood of $\leq$ in the space of monomial orderings.
\end{proof}

As an application of Theorem \ref{compact} we can give now a very 
short proof of the existence of universal Gr\"obner bases 
(following the idea of the proof given in \cite{Sc}).

\begin{theorem}\cite{Sc} 
For any ideal $I\triangleleft \K[x_1,...,x_n]$ there exists a finite set of
polynomials $f_1,...,f_s\in I$ which is a Gr\"obner basis of $I$
with respect to any monomial ordering. Such set is called
a universal Gr\"obner basis.
\end{theorem}

\begin{proof} 
For any $f_1,...,f_s\in I$ let $V_{\{f_1,...,f_s\}}$ be the set of
monomial orderings for which $\{f_1,...,f_s\}$ is a Gr\"obner basis
of $I.$ By Proposition \ref{open}, $V_{\{f_1,...,f_s\}}$ is open.
Since, by Corollary \ref{well_orders_closed}, 
the space of monomial orderings is compact, it has a finite cover by 
$V_{\{f_1,...,f_s\}},...,V_{\{h_1,...,h_t\}}.$
The set $\{f_1,...,f_s,...,h_1,...,h_t\}$ is the universal 
Gr\"obner basis of $I.$
\end{proof}

The following example shows that Corollary \ref{well_orders_closed} 
(which was an essential part of the above proof) does not hold for 
some monoids other then $\Z_{\geq 0}^n.$ 

{\bf Example:} The monoid of monomials in $\K[x_1,x_2,....]$
is isomorphic to the infinite sum, $\Z_{\geq 0}^\infty.$
Each linear order, $<,$ on the set $\{x_1,x_2,...\}$ induces a
lexicographic order on $\K[x_1,x_2,....]:$ Let
$m=x_1^{m_1}x_2^{m_2}...$ and $n=x_1^{n_1}x_2^{n_2}...$ be two arbitrary
monomials in $\K[x_1,x_2,....].$ Let $x_i$ be the smallest variable
(with respect to $<$) for which $m_i\ne n_i.$ Then
$m<n$ if and only if $m_i<n_i.$ We denote such induced order on 
$\K[x_1,x_2,....]$ by $<_L$ (L stands for lexicographic).
$<_L$ is a well-order if and only if $<$ is a well-order.
Consider now a sequence of well-orders $<^1,<^2,...$ such that
$$x_n<^n x_{n-1}<^n ... <^n x_2 <^n x_1 <^n x_{n+1} <^n x_{n+2} <^n x_{n+3}
 <^n ...$$
Consider a filtration $\{G_*\}$ of $\Z_{\geq 0}^\infty$ such that
$G_n$ contains all monomials in variables $x_1,...,x_n$ of total degree
at most $n.$ It is easy to check that the sequence of
lexicographic orders $<^1_L,<^2_L,...$ converges with respect to the
metric induced by the filtration $\{G_*\}$ to a new lexicographic order
$<^\infty_L$ induced by 
$$... <^\infty x_n <^\infty  x_{n-1}<^\infty ... <^\infty x_2 <^\infty x_1.$$
Hence $<^\infty_L$ is not a well-order on $\K[x_1,x_2,....],$ and therefore,
$LWO(\Z_{\geq 0}^\infty)$ is not a closed subset of $LO(\Z_{\geq 0}^\infty).$

\section{Acknowledgments}
We thank R. Bryant for suggesting a simplification of the proof 
of Proposition \ref{standard}(2) and the referee of this paper for
many insightful comments.

\vspace*{1in}

\centerline{Dept. of Mathematics, 244 Math Bldg, SUNY Buffalo, NY 14260} 
\centerline{asikora@buffalo.edu}

\end{document}